\documentclass [12pt]{article}
\usepackage{amsfonts}
\usepackage{amssymb}
\usepackage{multirow}
\usepackage{amsmath,amssymb,amsthm}

\newtheorem{remark}{Remark}
\newtheorem{lemma}{Lemma}
\newtheorem{proposition}{Proposition}
\newtheorem{corollary}{Corollary}
\newcommand{\re}[1]{\mbox{${\rm (\ref{#1})}$}}

\newcommand{\mysection}[1]{\section{#1}\setcounter{equation}{0}}

\title{Optimal upper bounds on expected $k$th record values\\ from IGFR distributions}
\author{Agnieszka Goroncy\\
Nicolaus Copernicus University\\
Chopina 12/18, 87100 Toru\'n, Poland\\
e-mail: gemini@mat.umk.pl
}

\date{}
\begin{document}

\maketitle{}
\begin{abstract}
The paper concerns the optimal upper bounds on the expectations of the $k$th record values ($k\geq 1$) centered about the sample mean. We consider the case, when the records are based on the infinite sequence of the independent identically distributed random variables, which distribution function is restricted to the family of distributions with the increasing generalized failure rate (IGFR). Such a class can be defined in terms of the convex orders of some distribution functions. Particularly important examples of IGFR class are the distributions with the increasing density (ID) and increasing failure rate (IFR). Presented bounds were obtained with use of the projection method, and are expressed in the scale units based on the standard deviation of the underlying distribution function.\\


\vspace{2ex}

\end{abstract}

\mysection{Introduction}
Let us consider the infinite sequence $\{X_n\}$, $n\geq 1$, of independent and identically distributed random variables with the common cummulative distribution function (cdf) $F$ and finite mean $\mu$ and variance $\sigma^2$. By $X_{1:n},\ldots,X_{n:n}$ denote the order statistics of $X_1,\ldots,X_n$. Further, we are interested in the increasing subsequences of $X_1,X_2,\ldots$ of the $k$th greatest order statistics, for a fixed $k=1,2,\ldots$. Formally, we define the (upper) $k$th records $R^{(k)}_n$, $n=0,1,\ldots$, by introducing first the $k$th record times $T_n^{(k)}$ as
\begin{equation}\nonumber
T_0^{(k)}=k, \quad T_n^{(k)}=\min\{j>T_{n-1}^{(k)}: X_j>X_{T_{n-1}^{(k)}+1-k:T_{n-1}^{(k)}}\}, \quad n=1,2,\ldots.
\end{equation}
Then the $k$th record values are given by
\begin{equation}\nonumber
R_n^{(k)}=X_{T_n^{(k)}+1-k:T_n^{(k)}}, \quad n=0,1,\ldots.
\end{equation}
Note that classic upper records are defined by $k=1$, and we say that such a record occurs at time $j$ if $X_j$ is greater than the maximum of previous observations $X_1,\ldots,X_{j-1}$.

Records are widely used, not only in the statistical applications. The most obvious one that arises at the first glance is the prediction of sport achievements and natural disasters. The first mention of the classic records comes from Chandler (1952), while the $k$th record values were introduced by Dziubdziela and Kopoci\'nski (1976). For the comprehensive overview of the results on the record values the reader is referred to Arnold, Balakrishnan and Nagaraja (1998) and Nevzorov (2001).

The distribution function of the $k$th record value $R_n^{(k)}$ is given by the following formula
\begin{equation}\label{F_nk}
F_n^{(k)}=1-[1-F(x)]^k\sum\limits_{i=0}^{n}\frac{k^i}{i!}(-\ln[1-F(x)])^i.
\end{equation}
If the cdf $F$ is absolutely continuous with a probability density function (pdf) $f$, then
the distribution function (\ref{F_nk}) also has the pdf given by
$$
f_n^{(k)}(x)=\frac{k^{n+1}}{n!}(-\ln[1-F(x)])^{n} [1-F(x)]^{k-1}
f(x).
$$
In particular case of the standard unform underlying cdf $F$, the corresponding distribution of the uniform $k$th record is given by the following equations
\begin{eqnarray}
G^{(k)}_n(x)&=&1-(1-x)^k
\sum\limits_{i=0}^{n}\frac{k^i}{i!}[-\ln(1-x)]^i,\quad 0<x<1,
\quad k \geq 1,\; n \geq 0,\nonumber\\
g^{(k)}_n(x)&=&\frac{k^{n+1}}{n!}[-\ln(1-x)]^{n}(1-x)^{k-1},\quad
0<x<1, \quad k \geq 1,\; n \geq 0.\label{g_nk}
\end{eqnarray}

Now, recall the cdf of the generalized Pareto distribution
\begin{equation}\label{W_alpha}
W_\alpha(x)=\left\{
\begin{array}{ll}
1-(1-\alpha x)^{1/\alpha}, & \textrm{for } x\geq 0, \textrm{ if } \alpha<0,\\
&\\
1-(1-\alpha x)^{1/\alpha},&  \textrm{for } 0\leq x\leq \frac{1}{\alpha}, \textrm{ if } \alpha>0, \\
&\\
1-\textrm{e}^{-x}, & \textrm{for } x\geq 0, \textrm{ if } \alpha=0.
\end{array}
\right.
\end{equation}
Next, we say that the cdf $F$ precedes the cdf $G$ in the convex transform order, and we write $F\prec G$, if the composition $F^{-1}\circ G$ is concave on the support of $G$. Following the reasoning of Goroncy and Rychlik (2015) and Bieniek and Szpak (2017), we consider the following family of distributions with the increasing generalized failure rate defined as with respect to \re{W_alpha}
\begin{equation}
\textrm{IGFR}(\alpha)=\{F: F\prec W_\alpha\}.
\end{equation}
Indeed, if the distribution function $F$ is continuous with the density function $f$, then the generalized failure rate defined as
\begin{equation}\label{gfr}
\gamma_\alpha(x)=(W^{-1}_\alpha(F))'(x)=(1-F(x))^{\alpha-1}f(x),
\end{equation}
is increasing. Note that the expression in \re{gfr} is just the product of the conventional failure rate and a power of the survival function $1-F$.

For $\alpha=1$ we obtain the standard uniform distribution function $W_1=U$ and the family IGFR(1)=ID of the increasing density distributions, respectively. On the other hand for $\alpha=0$, the cdf $W_0=V$ is the cdf of the standard exponential distribution and in the result we get the family IGFR(0)=IFR of the increasing failure rate distributions.

The aim of this paper is to establish the optimal upper bounds on
\begin{equation}\label{problem}
\mathbb{E}\frac{R_n^{(k)}-\mu}{\sigma},
\end{equation}
where the cdf $F$ is restricted to the IGFR($\alpha$) class of distributions, for arbitrarily chosen $k\geq 1$ and $n\geq 1$. In the special case $n=0$, which reduces to the order statistics $X_{1:k}$, the readers are referred to Rychlik (2014), who established the optimal bounds for ID and IFR distributions.

The bounds on the $k$th records, in particular classic record values have been widely considered in the literature, beginning with Nagaraja (1978). He used the Schwarz inequality to obtain the upper bounds on the expectations of the classic records, which were expressed in terms of the mean and standard deviation of the underlying distribution. The H\"older inequality was used by Raqab (2000), who presented more general bounds expressed in the scale units generated by the $p$th central absolute moments $\sigma_p$, $p\geq 1$. Also, he considered records from the symmetric populations. Differences of the consecutive record values (called record spacings) based on the general populations and from distributions with the increasing density and increasing failure rate were considered by Rychlik (1997). His results were generalized by Danielak (2005) into the arbitrary record increments.

The $k$th record values were considered by Grudzie\'n and Szynal (1985), who by use of the Schwarz inequality obtained non-sharp upper bounds expressed in terms of the population mean and standard deviation. Respective optimal bounds were derived by Raqab (1997), who applied the Moriguti (1953) approach. Further, the H\"older along with the Moriguti inequality were used by Raqab and Rychlik (2002) in order to get more general bounds. Gajek and Okolewski (2003) dealt with the expected $k$th record values based on the non-negative decreasing density and decreasing failure rate populations evaluated in terms of the population second raw moments. Results for the adjacent and non-adjacent $k$th records were obtained by Raqab (2004) and Danielak and Raqab (2004a). Evaluations for the se\-cond records from the symmetric populations were considered by Raqab and Rychlik (2004). Danielak and Raqab (2004b) presented the mean-variance bounds on the expectations of $k$th record spacings from the decreasing density and decreasing failure rate families of distributions. Further, Raqab (2007) considered second record increments from decreasing density families. Bounds for the $k$th records from decreasing generalized failure rate populations were evaluated by Bieniek (2007). Expected $k$th record values, as well as their differences from bounded populations were determined by Klimczak (2007), who expressed the bounds in terms of the lengths of the support intervals.

Regarding the lower bounds on the record values, there are not many papers concerning the problem, in opposite to the literature on the lower bounds for the order statistics and their linear combinations (see e.g. Goroncy and Rychlik (2006a), Goroncy and Rychlik (2006b, 2008), Rychlik (2007), Goroncy (2009)).
The lower bounds on the expected $k$th record values expressed in units generated by the central absolute moments of various orders, in the general case of the arbitrary parent distributions were presented by Goroncy and Rychlik (2011). There are also a few papers concerning the lower bounds on records indirectly, namely in the more general case of the generalized order statistics (Goroncy (2014), Bieniek and Goroncy (2017)).

Below we present a procedure which provides the basis of obtaining the optimal upper bounds on \re{problem} in the case of our interest. It is well known that
$$
\mathbb{E}R_n^{(k)}=\int_0^1F^{-1}(x)g^{(k)}_n(x)dx=\int_0^1
F^{-1}(x)\frac{k^{n+1}}{n!}[-\ln(1-x)]^{n}(1-x)^{k-1}dx,
$$
therefore
\begin{equation}\label{E}
\mathbb{E}\frac{R_n^{(k)}-\mu}{\sigma}=
\int_0^1\frac{F^{-1}(x)-\mu}{\sigma}[g^{(k)}_n(x)-1]dx.
\end{equation}
Due to the further application, we subtract 1 from $g^{(k)}_n(x)$ in the formula above, but one could replace it with an arbitrary constant.
Changing the variables in \re{E}, for a fixed, absolutely continuous cdf $W$ with the pdf $w$ on the support $[0,d)$, $d\leq\infty$, we obtain
\begin{equation}\label{problem1}
\mathbb{E}\frac{R_n^{(k)}-\mu}{\sigma}=
\int_0^d\frac{F^{-1}(W(x))-\mu}{\sigma}(g^{(k)}_n(W(x))-1) w(x)dx.
\end{equation}
Further assume that $W$ satisfies
\begin{equation}\label{warunek_x2}
\int\limits_0^d x^2 w(x)dx<\infty.
\end{equation}
Let us consider the Hilbert space $L^2_W$ of the square integrable functions with respect to $w$ on $[0,d)$, and denote the norm of an arbitrary function $f\in L^2_W$ as
$$
||f||_W=\left(\int_0^d|f(x)|^2w(x)dx\right)^{1/2}.
$$
Moreover, let $P_W$ stand for the projection operator onto the following convex cone
\begin{equation}\label{C_W}
\mathcal{C}_W=\{g\in L^2_W: \textrm{$g$ is nondecreasing and concave}\}.
\end{equation}
In order to find the upper optimal bounds on \re{problem1}, we will use the Schwarz inequality combined with the well-known projection method (see Rychlik (2001), for details). It is clear that \re{problem1} can be bounded by the $L^2_{W}$-norm of the projection $P_{W}$ of the function $h_W(x)=g^{(k)}_n(W(x))-1$, as follows
\begin{equation}\label{problem_Ph}
\mathbb{E}\frac{R_n^{(k)}-\mu}{\sigma}\leq ||P_W h_W||_{W},
\end{equation}
with the equality attained for cdf $F$ satisfying
\begin{equation}\label{warunek_rownosci}
\frac{F^{-1}(W(x))-\mu}{\sigma}=\frac{P_W h_W(x)}{||P_W h_W||_{W}}.
\end{equation}
In our case we fix $W=W_\alpha$ and the problem of establishing the optimal upper bounds on \re{problem} easily boils down to determining the $L^2_{W_\alpha}$-norm of the projection $P_\alpha=P_{W_\alpha}$ of the function $h_{W_\alpha}$ onto $\mathcal{C}_{W_\alpha}$. Note that in order to apply the projection method, we need the condition \re{warunek_x2} to be fulfilled by the distribution function \re{W_alpha}. Bieniek (2008) showed, that in that case we need to confine ourselves to parameters $\alpha>-\frac{1}{2}$, what we do in our further considerations.

\mysection{Auxiliary results}
In this section we recall the results of Goroncy and Rychlik (2015, 2016), who determined the projection of the function $h\in L^2_W$ satisfying particular conditions, onto the cone \re{C_W} of nondecreasing and concave functions. These conditions are presented below.

(A) Let $h$ be bounded, twice differentiable function on $[0,d)$, such that
\begin{equation}\nonumber
\int\limits_0^d h(x)w(x)dx=0.
\end{equation}
Moreover, assume that $h$ is strictly decreasing on $(0,a)$, strictly convex increasing on $(a,b)$, strictly concave increasing on $(b,c)$ with $h(0)\leq 0<h(c)$, and strictly decreasing on $(c,d)$ with $h(d)=h(0)$ for some $0\leq a<b<c<d$.

The projection of the function $h$ satisfying conditions (A) onto the convex cone \re{C_W} is either first linear, then coinciding with $h$ and ultimately constant, or just linear and then constant, depending on the behaviour on some particular auxiliary functions, which are introduced below.

First, denote
\begin{equation}\label{T}
T_W(\beta)=h(\beta)(1-W(\beta))-\int\limits_\beta^dh(x)w(x)dx,\quad 0\leq \beta\leq d,
\end{equation}
which is decreasing on $(0,a)$, increasing on $(a,c)$ and decreasing on $(c,d)$, having the unique zero $\beta_*$ in $(a,c)$. Moreover, let
\begin{eqnarray}
\lambda_W(y)&=&\frac{\int\limits_0^y(x-y)(h(x)-h(y))w(x)dx}{\int\limits_0^y(x-y)^2w(x)dx},\label{lambda}\\
Y_W(y)&=& \lambda_W(y)-h'(y),\label{Y}\\
Z_W(y)&=& \int\limits_0^y(h(x)-h(y)-\lambda_W(y)(x-y))w(x)dx,\label{Z}
\end{eqnarray}
for $0\leq y<d$. The precise form of the projection of the function $h$ satisfying (A) onto the cone $\re{C_W}$ is described in the proposition below (cf. Goroncy and Rychlik (2016), Proposition 1).
\begin{proposition} \label{Proposition_GR}
If the zero $\beta_*\in(a,c)$ of \re{T} belongs to the interval $(b,c)$ and the set $\mathcal{Y}=\{y\in(b,\beta_*): Y_W(y)\geq 0, Z_W(y)=0\}$ is nonempty, then
\begin{equation}\nonumber
P_Wh(x)=\left\{
\begin{array}{lll}
h(y_*)+\lambda_W(y_*)(x-y_*),&0\leq x<y_*,\\
h(x),&y_*\leq x<\beta_*,\\
h(\beta_*), &\beta_*\leq x<d,
\end{array}
\right.
\end{equation}
where $y_*=\inf\{y\in \mathcal{Y}\}$ is the projection of h onto \re{C_W}. Otherwise we define
\begin{equation}\nonumber
  P_yh(x)=\frac{\int\limits_y^dh(t)w(t)dt}{1-W(y)}\left[\frac{(x-y)\mathbf{1}_{(0,y)}(x)}{-\int\limits_0^y(t-y)w(t)dt}+1\right], \quad \beta_*\leq y<d,
\end{equation}
with
\begin{equation}\nonumber
||P_yh||^2_W=\frac{\left(\int\limits_y^d h(x)w(x)dx\right)^2\left[\int\limits_0^y(x-y)^2w(x)dx-\left(\int\limits_0^y(x-y)w(x)dx)\right)^2\right]}{\left(\left(1-W(y)\right)\int\limits_0^y(x-y)w(x)dx\right)^2}.
\end{equation}
Let $\mathcal{Z}$ denote the set of arguments $y\geq \beta_*$ satisfying the following condition
\begin{equation}\label{warunek_na_y}
\frac{\int\limits_y^dh(x)w(x)dx}{1-W(y)}=-\frac{\int\limits_0^y(x-y)h(x)w(x)dx\int\limits_0^y(x-y)w(x)dx}{\int\limits_0^y(x-y)^2w(x)dx-\left(\int\limits_0^y(x-y)w(x)dx\right)^2}>0.
\end{equation}
Then $\mathcal{Z}$ is nonempty and $P_Wh(x)=P_{y_*}h(x)$ for unique $y_*=\arg\max\limits_{y\in\mathcal{Z}}||P_yh||^2_W$.
\end{proposition}
Note that there are only two possible shapes of projection functions of the function $h$ onto \re{C_W}. The first one requires compliance with certain  conditions and can be briefly described as: linear - identical with $h$ - constant (l-h-c, for short). The second possible shape does not have a part which corresponds to the function $h$, and we will refer to it as l-c (linear and constant) from now on.
The original version of this proposition can be found in Goroncy and Rychlik (2015), however there was no clarification about the parameter $y$ in case of the l-c type of the projection, therefore we refer to Goroncy and Rychlik (2016).

We will also need some results on the projection of the functions satisfying conditions ($\tilde{A}$), which are a slight modification of conditions (A). We state that the function $h$ satisfies ($\tilde{A}$) if conditions (A) are modified so that $c=d=\infty$ and $\sup\limits_{x>0} h(x)=\lim\limits_{x\rightarrow \infty}h(x)>0$. This in general means that the function does not have the decreasing part at the right end of the support and in particular does not have to be bounded from above. The proposition below (cf. Goroncy and Rychlik (2016), Proposition 6) describes the shape of the projection in this case.

\begin{proposition}\label{Proposition_GR2}
If the function $h$ satisfies conditions $(\tilde{A})$, then the set $\mathcal{\tilde{Y}}=\{y>b: Y_W(y)\geq 0, Z_W(y)=0\}$ is nonempty and for $y_*=\inf\{y\in\mathcal{\tilde{Y}}\}$ we have
\begin{equation}\nonumber
  P_Wh(x)=\left\{
  \begin{array}{ll}
  h(y_*)+\lambda_W(y_*)(x-y_*),& 0\leq x<y_*,\\
  h(x),& x\geq y_*.
  \end{array}
  \right.
\end{equation}

\end{proposition}

\mysection{Main results}
Let us focus now on the case $W=W_\alpha$ and denote
\begin{equation}\label{h_alpha}
h_\alpha(x)=h_{W_\alpha}(x)=\hat{g}^{(k)}_n(x)-1,
\end{equation}
where
\begin{equation}\label{g_hat}
\hat{g}^{(k)}_n=g^{(k)}_n\circ W_\alpha.
\end{equation}
We also denote $\hat{G}^{(k)}_n(x)=G^{(k)}_n(W_\alpha(x))$.

The substantial matter in determining the bounds on \re{problem} is to learn the shapes of the functions \re{h_alpha} for arbitrary $k=1,2,\ldots$ and $n\geq 1$, which correspond with the shapes of compositions \re{g_hat}, and are presented in the lemma below (comp. with Bieniek (2007), Lemma 3.2).

\begin{lemma} \label{lemma_shape}
If $\alpha>0$, then the shape of $\hat{g}^{(k)}_n$ is as follows:
\begin{itemize}
\item[(i)] If $k=1$, then $\hat{g}^{(1)}_n$, $n\geq 1$, is convex increasing.
\item[(ii)] If $2\leq k<\frac{\alpha}{2}+1$, then $\hat{g}^{(k)}_n$, $n\geq 1$, is convex increasing, concave increasing and concave decreasing.
\item[(iii)] If $\frac{\alpha}{2}+1\leq k\leq 1+\alpha$, then $\hat{g}^{(k)}_1$ is concave increasing-decreasing, and $\hat{g}^{(k)}_n$, $n\geq 2$, is convex increasing, concave increasing and concave decreasing.
\item[(iv)] If $k>1+\alpha$, then $\hat{g}^{(k)}_1$ is concave increasing, concave decreasing and convex decreasing, and $\hat{g}^{(k)}_n$, $n\geq 2$, is convex increasing, concave increasing, concave decreasing and convex decreasing.
\end{itemize}
If $\alpha=0$, then the shape of $\hat{g}^{(k)}_n$ is as follows:
\begin{itemize}
\item[(i)] If $k=1$, then $\hat{g}^{(1)}_1$ is linear increasing and $\hat{g}^{(1)}_n$, $n\geq 2$, is convex increasing.
\item[(ii)] If $k\geq 2$, then $\hat{g}^{(k)}_1$ is concave increasing and then decreasing, $\hat{g}^{(k)}_n$, $n\geq 2$, is convex increasing, concave increasing, and decreasing.
\end{itemize}
If $\alpha<0$, then the shape of $\hat{g}^{(k)}_n$ is as follows:
\begin{itemize}
  \item[(i)] If $k=1$, then $\hat{g}^{(1)}_1$ is concave increasing, $\hat{g}^{(1)}_n$, $n\geq 2$, is convex increasing and concave increasing.
  \item[(ii)] If $k\geq 2$, then $\hat{g}^{(k)}_1$ is concave increasing, concave decreasing and convex decreasing and $\hat{g}^{(k)}_n$, $n\geq 2$, is convex increasing, concave increasing, concave decreasing and convex decreasing.
\end{itemize}
\end{lemma}
It is worth mentioning that slight differences between the lemma above and Lemma 3.2 in Bieniek (2007) are the result of different notations of record values.

Note that the case $\alpha=1$ is covered by the above lemma, except the setting (ii) for $\alpha>0$, which is not possible in this case (cf. Rychlik (2001), p.136). Case $\alpha=0$ comes from Rychlik (2001, p.136). In order to determine the shape of \re{g_hat} for $\alpha\neq 0$, we notice that for $k\geq 2$ and $n\geq 2$
\begin{eqnarray}
(\hat{g}^{(k)}_n(x))'&=& \frac{1}{1-\alpha x}[k\hat{g}^{(k)}_{n-1}(x)-(k-1)\hat{g}^{(k)}_n(x)],\\
(\hat{g}^{(k)}_n(x))''&=&\frac{1}{(1-\alpha x)^2}[k^2\hat{g}^{(k)}_{n-2}(x)-k(2k-\alpha-2)\hat{g}^{(k)}_{n-1}(x)\nonumber\\
&&+(k-1)(k-1-\alpha)\hat{g}^{(k)}_n(x)],
\end{eqnarray}
and use the variation diminishing property (VDP) of the linear combinations of $g^{(k)}_n\circ W_\alpha$ (see Gajek and Okolewski, 2003). Other special cases of $k, n$ we calculate separately in order to obtain the shapes of \re{g_hat}.

Faced with this knowledge, we conclude that \re{h_alpha} satisfy conditions (A) with $a=0$, for $k\geq 2$ and $n\geq 2$ if $\alpha=0$, for $n=1$ and $2\leq k<\frac{\alpha}{2}+1$ or $k\geq 2$ and $n\geq 2$ if $\alpha>0$, as well as for $k\geq 2$ and $n\geq 2$ if $\alpha<0$. Moreover, we have
$$
\int\limits_0^d h_\alpha(x)w(x)dx=\int\limits_0^1 (g^{(k)}_n(u)-1)du=0,
$$
$h_\alpha(0)=g^{(k)}_n(W_\alpha(0))-1=-1$, $h_\alpha(d)=g^{(k)}_n(W_\alpha(d))-1=-1$. The value of \re{h_alpha} in the local maximum point $c$ has to be positive, since the function starts and finishes with negative values and integrates to zero, which means that \re{h_alpha} has to cross the $x$-asis and changes the sign from negative to positive, finishing with negative value at $d$.
Therefore, we can use Proposition \ref{Proposition_GR} in order to obtain the projection of \re{h_alpha} onto \re{C_W} and finally determine the desired bounds according to \re{problem_Ph}. Moreover, \re{h_alpha} satisfy conditions ($\tilde{A}$) in case of the first record values ($k=1$) for $\alpha<0$, $n\geq 2$ and we are entitled to use Proposition \ref{Proposition_GR2} then. Other cases can be dealt without the above results. These imply the particular shapes of the projections which can be one of the three possible kinds. The first one coincides with the original function \re{h_alpha} (first values of the classic records for $\alpha\leq 0$), the second shape is the linear increasing function (classic record values for $n\geq 1$ and $\alpha>0$ or $n\geq 2$ and $\alpha=0$), and the last one is the projection coinciding with the function $h$ at the beginning and ultimately constant (first values of the $k$th records for $\alpha<0$ and $k\geq 2$ or $\alpha>0$, $k\geq \max\{2,\frac{\alpha}{2}+1\}$).

In order to simplify the notations, we will denote the projection of function \re{h_alpha} onto \re{C_W} with respect to \re{W_alpha} by $P_\alpha h_\alpha$ from now on.\\

\subsection{Bounds for the classic records}
In the proposition below we present the bounds on the classic record values ($k=1$). This case does not require using the Proposition \ref{Proposition_GR}, since the shapes of the densities of records do not satisfy conditions (A), but possibly satisfy conditions ($\tilde{A}$).
\begin{proposition}\label{Proposition_k1} Assume that $k=1$.
\begin{itemize}
\item [(i)] Let $-\frac{1}{2}<\alpha<0$. If $n=1$, then we have the following bound
\begin{equation}\label{oszac_rowne_h}
 \mathbb{E}\frac{R_1^{(1)}-\mu}{\sigma}\leq 1,
\end{equation}
with the equality attained for the exponential distribution function
\begin{equation}\label{rozklad_rownosc_n1}
F(x)=1-\exp\left\{-1-\frac{x-\mu}{\sigma}\right\}, \quad x> \mu-\sigma.
\end{equation}
If $n\geq 2$, then the set $\tilde{\mathcal{Y}}=\{y>b: Y_\alpha(y)\geq 0, Z_\alpha(y)=0\}$ is nonempty and for $y_*=\inf\{y\in\tilde{\mathcal{Y}}\}$ we have
\begin{equation}\label{oszac_ogolne3}
 \mathbb{E}\frac{R_n^{(1)}-\mu}{\sigma}\leq C_{\alpha}(y_*),
\end{equation}
where
\begin{eqnarray}
C^2_{\alpha}(y)&=&\left[1-(1-\alpha y)^{1/\alpha}\right]\left[1+\left(\hat{g}^{(1)}_n(y)-1\right)^2\right]-\hat{G}^{(1)}_n(y) \nonumber\\
&&\frac{2(1+2\alpha)\left(\hat{g}^{(1)}_n(y)-1\right)[1-y(1+\alpha)-(1-\alpha y)^{1/\alpha+1}]}{2-2(1-\alpha y)^{1/\alpha+2}+y(1+2\alpha)(\alpha y+y-2)}\nonumber\\
&&\cdot \left[\hat{g}^{(1)}_n(y)\frac{y(1+\alpha)-1+(1-\alpha y)^{1/\alpha+1}}{1+\alpha}-\int\limits_0^y\hat{G}^{(1)}_n(x)dx\right]\nonumber\\
&& +\frac{}{}\frac{(1-2\alpha)\left([y(1+\alpha)-1+(1-\alpha y)^{1/\alpha+1}]\hat{g}^{(1)}_n(y)-(1+\alpha)\int\limits_0^y \hat{G}^{(1)}_n(x)dx\right)^2}{(1+\alpha)[2-2(1-\alpha y)^{1/\alpha+2}+y(1+2\alpha)(\alpha y+y-2)]}.\nonumber\\
&&\label{C_kwadrat}
\end{eqnarray}
The equality in \re{oszac_ogolne3} is attained for distribution functions $F\in$ IGFR$(\alpha)$ that satisfy the following condition
\begin{equation}\nonumber
F^{-1}(W_\alpha(x))=\left\{
\begin{array}{ll}
\dfrac{\sigma}{C_{\alpha}(y_*)}\left[\hat{g}^{(1)}_n(y_*)-1+\lambda_\alpha(y_*)(x-y_*)\right]+\mu,&0\leq x<y_*,\\
&\\
\dfrac{\sigma}{C_{\alpha}(y_*)}\left[\hat{g}^{(1)}_n(x)-1\right]+\mu,&y_*\leq x<d.\\
\end{array}
\right.
\end{equation}
\item[(ii)] Let now $\alpha=0$. We have the following bound
\begin{equation}\nonumber
 \mathbb{E}\frac{R_n^{(1)}-\mu}{\sigma}\leq n,
\end{equation}
with the equality attained for the exponential distribution function \re{rozklad_rownosc_n1}.
\item [(iii)] Suppose $\alpha> 0$, $n\geq 1$. Then we have the following bound
\begin{equation}\label{oszac_liniowe}
 \mathbb{E}\frac{R_n^{(1)}-\mu}{\sigma}\leq \sqrt{\dfrac{(2\alpha+1)(2a_*b_*+(\alpha+1)b_*^2)+2a_*^2}{(1+\alpha)(2\alpha+1)}},
\end{equation}
where
\begin{eqnarray}
a_*&=&(1+\alpha)^2(2\alpha+1)\left[\frac{1}{\alpha(1+\alpha)}-\int\limits_0^{1/\alpha}\hat{G}^{(1)}_n(x)dx\right],\label{a_*}\\
b_*&=&-\frac{a_*}{1+\alpha}=-(1+\alpha)(2\alpha+1)\left[\frac{1}{\alpha(1+\alpha)}-\int\limits_0^{1/\alpha}\hat{G}^{(1)}_n(x)dx\right].\label{b_*}
\end{eqnarray}
The equality in \re{oszac_liniowe} is attained for the following distribution function
$$
F(x)=1-\left(1-\alpha\left(1+\dfrac{(x-\mu)}{\sigma a_*}\sqrt{\dfrac{(2\alpha+1)(2a_*b_*+(\alpha+1)b_*^2)+2a_*^2}{(1+\alpha)(2\alpha+1)}}\right)\right)^{1/\alpha}.
$$
\end{itemize}
\end{proposition}

\textsc{Proof.} Fix $k=1$.
 Let us first consider case (i), i.e. $\alpha<0$. Here we have to add an additional restriction $\alpha>-\frac{1}{2}$, which has been mentioned at the end of Section 2.
If $n=1$, then the function $h_\alpha$ is increasing and concave, hence its projection onto \re{C_W} is the same as $h_\alpha$. The bound can be determined via its norm, which square is given by
\begin{equation}\label{3.25}
||P_\alpha h_\alpha||^2=||h_\alpha||^2=\int\limits_0^d(\hat{g}^{(k)}_n(x)-1)^2w_\alpha(x)dx=\frac{k^{2(n+1)}(2n)!}{(n!)^2(2k-1)^{2n+1}}-1,
\end{equation}
since
\begin{eqnarray}
  \int\limits_{0}^d\left[\hat{g}^{(k)}_n(x)\right]^2w_\alpha(x)dx&=&\int\limits_{0}^1\left[g^{(k)}_n(x)\right]^2dx=\nonumber\\
  &=&\frac{k^{2(n+1)}(2n)!}{(n!)^2(2k-1)^{2n+1}}\int\limits_0^1g^{(2k-1)}_{2n}(x)dx=\frac{k^{2(n+1)}(2n)!}{(n!)^2(2k-1)^{2n+1}}.\nonumber
\end{eqnarray}

Taking into account that $k$ as well as $n$ are equal to one, formula \re{3.25} implies \re{oszac_rowne_h}.\\
Suppose now that $n\geq 2$. Note that in this case \re{h_alpha} satisfy conditions ($\tilde{A}$). Using Proposition \ref{Proposition_GR2}, we have the following projection of \re{h_alpha} onto the cone \re{C_W},
\begin{equation}\nonumber
P_\alpha h_\alpha(x)=\left\{
\begin{array}{ll}
\hat{g}^{(k)}_n(y_*)-1+\lambda_\alpha(y_*)(x-y_*),&0\leq x<y_*,\\
&\\
\hat{g}^{(k)}_n(x)-1,& y_*\leq x<d.
\end{array}
\right.
\end{equation}
An appropriate counterpart of function \re{lambda} in our case is
\begin{equation}
\lambda_\alpha(y)=\lambda_{W_\alpha}(y)=\frac{\hat{g}^{(k)}_n(y)\int\limits_0^yW_\alpha(x)dx-\int\limits_0^y\hat{G}^{(k)}_n(y)dx}{2y\int\limits_0^yW_\alpha(x)dx-2\int\limits_0^yxW_\alpha(x)dx},\label{lambda_alpha}\\
\end{equation}
with $k=1$, since simple calculations show that
\begin{eqnarray}
  \int\limits_0^y(x-y)w_\alpha(x)dx &=& -\int\limits_0^yW_\alpha(x)dx, \label{calka_x-y}\\
  \int\limits_0^y(x-y)^2w_\alpha(x)dx &=&  2\left(y\int\limits_0^yW_\alpha(x)dx-\int\limits_0^yxW_\alpha(x)dx\right),\label{calka_x-y_kw}\\
  \int\limits_0^yx\hat{g}^{(k)}_n(x)w_\alpha(x)dx&=&y\hat{G}^{(k)}_n(y)-\int\limits_0^y\hat{G}^{(k)}_n(x)dx.\label{calka_g_hat_w}
\end{eqnarray}
Having
\begin{eqnarray}
\int\limits_0^y W_\alpha(x)dx&=&y+\frac{(1-\alpha y)^{1+1/\alpha}-1}{1+\alpha},\nonumber\\
\int\limits_0^y xW_\alpha(x)dx&=&\frac{1}{2}y^2+\frac{(1-\alpha y)^{1+1/\alpha}}{1+\alpha}y+\frac{(1-\alpha y)^{2+1/\alpha}-1}{(1+\alpha)(1+2\alpha)},\nonumber
\end{eqnarray}
for $\alpha\neq 0$, we conclude that \re{lambda_alpha} takes the form
\begin{equation}\label{lambda_alpha_ost}
\lambda_\alpha(y)=(1+2\alpha)\frac{[y(1+\alpha)-1+(1-\alpha y)^{1/\alpha+1}]\hat{g}^{(k)}_n(y)-(1+\alpha)\int\limits_0^y \hat{G}^{(k)}_n(x)dx}{2-2(1-\alpha y)^{1/\alpha+2}+y(1+2\alpha)(\alpha y+y-2)},
\end{equation}
with $k=1$. In consequence for $-\frac{1}{2}<\alpha<0$ and $k=1$ we have
$$
||P_\alpha h_\alpha||^2=C^2_\alpha(y_*),
$$
where $C^2_\alpha(y)$ is given in \re{C_kwadrat}. The square root of the expression above determines the optimal bound on \re{problem}.\\

Consider now case (iii) with $\alpha> 0$ and $n\geq 1$, which requires more explanation. With such parameters function $h_\alpha$ is increasing and convex. This implies that its projection onto the cone of the nondecreasing and concave functions is linear increasing. The justification for this is similar as e.g. in Rychlik (2014, p.9). The only possible shape of the closest increasing and convex function to the function $h_\alpha$ is the linear increasing one $P_\alpha h_\alpha(x)=a_0x+b_0$, say, which has at most two crossing points with $h_\alpha$. Since
\begin{equation}\label{warunek_calka_rzutu}
\int\limits_0^dh_\alpha(x)w_\alpha(x)dx=\int\limits_0^dP_\alpha h_\alpha(x)w_\alpha(x)=0,
\end{equation}
(see e.g. Rychlik (2001)), we obtain
\begin{equation}\label{b_0}
b_0=-\dfrac{a_0}{1+\alpha}.
\end{equation}
Next, in order to determine the optimal parameter $a_0$, we need to minimize the distance between the function \re{h_alpha} and its projection
\begin{equation}\nonumber
D_\alpha(a_0)=||P_\alpha h_\alpha-h_\alpha||^2.
\end{equation}
For $\alpha>0$ we have $\hat{g}^{(1)}_n=\frac{1}{n!}(-\ln(1-\alpha x)^{1/\alpha})^n$ and $w_\alpha(x)=(1-\alpha x)^{1/\alpha-1}$. Therefore
\begin{equation}\label{D_alpha}
D_\alpha(a)=a^2\int\limits_0^{1/\alpha}\left(x-\frac{1}{1+\alpha}\right)^2w_\alpha(x)dx-2a\int\limits_{0}^{1/\alpha}\left(x-\frac{1}{1+\alpha}\right)h_\alpha(x)w_\alpha(x)dx+\int\limits_0^{1/\alpha}h_\alpha^2(x)w_\alpha(x)dx.
\end{equation}
Using
\begin{eqnarray}
\int\limits_0^{1/\alpha}\left(x-\frac{1}{1+\alpha}\right)^2w_\alpha(x)dx&=&\frac{1}{(2\alpha+1)(1+\alpha)^2},\nonumber\\
\int\limits_{0}^{1/\alpha}\left(x-\frac{1}{1+\alpha}\right)h_\alpha(x)w_\alpha(x)dx&=&\frac{1}{\alpha(1+\alpha)}-\int\limits_0^{1/\alpha}\hat{G}^{(1)}_n(x)dx,\nonumber
\end{eqnarray}
we get the minimum of \re{D_alpha} equal to \re{a_*}. Since \re{b_0}, we also obtain \re{b_*}. Finally, the optimal bound can be determined by calculating the square root of
\begin{equation}\nonumber
||P_\alpha h_\alpha||^2=\int\limits_0^{1/\alpha}(a_*x+b_*)^2(1-\alpha x)^{1/\alpha-1}dx,
\end{equation}
which equals \re{oszac_liniowe}.\\

Let finally consider the case (ii). Note that for $n=1$ function $h_\alpha$ is increasing and concave, and the case is analogous to (i) with $n=1$, when we get the bound equal to 1.
If $n\geq 2$, then $h_\alpha$ is increasing and convex and its projection onto the cone of the nondecreasing and concave functions is linear, as in case (iii). Here the analogue to \re{b_0} is $b_0=-a_0$. For $\alpha=0$, we have $\hat{g}^{(1)}_n(x)=\frac{x^n}{n!}$, $n\geq 2$, which gives us the distance function $D_0(a)=a^2-2an-2+\frac{(2n)!}{(n!)^2}$, which is minimized for $a_0=n$. Hence $P_0h_0(x)=n(x-1)$, and we get the optimal bound equal to $||P_0h_0||=n$.\\

The distributions for which the equalities are attained in all the above cases can be determined using the condition \re{warunek_rownosci} with $W=W_\alpha$ and $h=h_\alpha$. \hfill$\Box$\\

\subsection{Bounds for the $k$th records, $k\geq 2$}
As soon as we give some auxiliary calculations, we are ready to formulate the results on the upper bounds of the expected $k$th record values, $k\geq 2$ based on the the IGFR($\alpha$) family of distributions.
For $W_\alpha$ being the GPD distribution, we have the corresponding function of \re{T} given by
\begin{equation}\label{T_alpha}
T_\alpha(\beta)=T_{W_\alpha}(\beta)=(1-W_\alpha(\beta))\left[-\frac{1}{k}\sum\limits_{i=0}^{n-1}\hat{g}^{(k)}_i(\beta)+\left(1-\frac{1}{k}\right)\hat{g}^{(k)}_n(\beta)\right], \beta\in(0,d).
\end{equation}
Knowing the properties of \re{T}, we notice that $T_\alpha(0)=-1$, $T_\alpha(d)=0$, and \re{T_alpha} is first negative, then positive, which follows from the VDP property of the linear combinations of $\hat{g}^{(k)}_i$, $i=0,\ldots,n$ (see Gajek and Okolewski, 2003). We conclude that \re{T_alpha} increases from $-1$ to $T_\alpha(c)>0$, and then decreases to $0$, which means that \re{T_alpha} has the unique zero $\beta_*$ in the interval $(0,c)$. Therefore the condition $\beta_*\in(b,c)$, required in the Proposition \ref{Proposition_GR} for the l-h-c type of the projection of the function \re{h_alpha}, is equivalent to $T_\alpha(b)<0$.\\
Moreover, respective functions \re{lambda}-\re{Z} in case of the $n$th values of the $k$th records, $n=1,\ldots$, $k=2,\ldots$, are following
\begin{eqnarray}
Y_\alpha(y)&=&Y_{W_\alpha}(y)=\frac{\hat{g}^{(k)}_n(y)\int\limits_0^yW_\alpha(x)dx-\int\limits_0^y\hat{G}^{(k)}_n(y)dx}{2y\int\limits_0^yW_\alpha(x)dx-2\int\limits_0^yxW_\alpha(x)dx}\nonumber\\
&&-w_\alpha(y)\frac{k^{n+1}}{(k-1)^n}\left(\hat{g}^{(k-1)}_{n-1}(y)-\hat{g}^{(k-1)}_n(y)\right),\nonumber\\
Z_\alpha(y)&=&Z_{W_\alpha}(y)=\hat{G}^{(k)}_n(y)-W_\alpha(y)\hat{g}^{(k)}_n(y)+\lambda_\alpha(y)\int\limits_0^yW_\alpha(x)dx,\nonumber
\end{eqnarray}
together with $\lambda_\alpha$ presented in \re{lambda_alpha}, since \re{calka_x-y}, \re{calka_x-y_kw} and \re{calka_g_hat_w} hold.

\begin{proposition}\label{Proposition_ogolne} Let $F\in$IGFR($\alpha$) and $k\geq 2$. \\
(i) Fix $n=1$ and let $-\frac{1}{2}<\alpha\leq 0$ or $0<\alpha\leq 2k-1$. Then for
\begin{equation}\label{beta_star}
\beta^*=\left\{
\begin{array}{ll}
\frac{1}{k(k-1)},&\alpha=0,\\
&\\
\frac{1}{\alpha}\left(1-{\rm exp}\left(-\frac{\alpha}{k(k-1)}\right)\right),&\alpha\neq0,
\end{array}
\right.
\end{equation}
we have the following bound
\begin{equation}\label{oszac_ogolne0}
\mathbb{E}\frac{R_1^{(k)}-\mu}{\sigma}\leq B_{\alpha,0}(\beta^*),
\end{equation}
where
\begin{equation}
B_{\alpha,0}^2(\beta)=\frac{2k^{4}}{(2k-1)^3}\hat{G}^{(2k-1)}_{2}(\beta)-2\hat{G}^{(k)}_1(\beta)+W_\alpha(\beta)+(1-W_\alpha(\beta))\left[\hat{g}^{(k)}_1(\beta)-1\right]^2.
\end{equation}
The equality in \re{oszac_ogolne0} is attained for the following IGFR$(\alpha)$ distribution function $F$
\begin{equation}
  F^{-1}(W_\alpha(x))=\left\{
  \begin{array}{ll}
  \dfrac{\sigma}{B_{\alpha,0}(\beta^*)}(\hat{g}^{(k)}_1(x)-1)+\mu,&0\leq x<\beta^*,\\
  &\\
  \dfrac{\sigma}{B_{\alpha,0}(\beta^*)}(\hat{g}^{(k)}_1(\beta_*)-1)+\mu,&\beta^*\leq x\leq d.
  \end{array}
  \right.
\end{equation}
(ii)
Let $-\frac{1}{2}<\alpha\leq 0$ and $n\geq 2$ or $\alpha>0$ and $n=1$ with $2\leq k<\frac{\alpha}{2}+1$ or $n\geq 2$.  \\
 If $T_\alpha(b)<0$, and the set $\mathcal{Y}_\alpha=\{y\in(b,\beta_*): Y_\alpha(y)\geq 0, Z_\alpha(y)=0\}$ is nonempty, then let $y_*=\inf\{y\in \mathcal{Y_\alpha}\}$ and we have the following bound
\begin{equation}\label{oszac_ogolne1}
\mathbb{E}\frac{R_n^{(k)}-\mu}{\sigma}\leq B_{\alpha,1}(y_*,\beta_*),
\end{equation}
where
\begin{eqnarray}
B_{\alpha,1}^2(y,\beta)&=&W_\alpha(y)\left[\hat{g}^{(k)}_n(y)-1\right]^2+(1-W_\alpha(\beta))\left[\hat{g}^{(k)}_n(\beta)-1\right]^2\nonumber\\
&&-2\lambda_\alpha(y)[\hat{g}^{(k)}_n(y)-1]\int\limits_0^yW_\alpha(x)dx+2\lambda_\alpha^2(y)\left[y\int\limits_0^yW_\alpha(x)dx-\int\limits_0^yxW_\alpha(x)dx\right]\nonumber\\
&&+\frac{k^{2(n+1)}(2n)!}{(n!)^2(2k-1)^{2n+1}}\left[\hat{G}^{(2k-1)}_{2n}(\beta)-\hat{G}^{(2k-1)}_{2n}(y)\right]\nonumber\\
&&-2\left[\hat{G}^{(k)}_n(\beta)-\hat{G}^{(k)}_n(y)\right]+W_\alpha(\beta)-W_\alpha(y),\nonumber
\end{eqnarray}
with $\lambda_\alpha$ given by \re{lambda_alpha}. The equality in \re{oszac_ogolne1} is attained for distributions $F\in$ IGFR($\alpha$) satisfying the following condition
\begin{equation}\nonumber
F^{-1}(W_\alpha(x))=\left\{
\begin{array}{ll}
\dfrac{\sigma}{B_{\alpha,1}(y_*,\beta_*)}\left[\hat{g}^{(k)}_n(y_*)-1+\lambda_\alpha(y_*)(x-y_*)\right] +\mu, & 0\leq x<y_*,\\
 &\\
\dfrac{\sigma}{B_{\alpha,1}(y_*,\beta_*)}\left(\hat{g}^{(k)}_n(x)-1\right) +\mu, & y_*\leq x<\beta_*,\\
&\\
\dfrac{\sigma}{B_{\alpha,1}(y_*,\beta_*)}\left(\hat{g}^{(k)}_n(\beta_*)-1\right) +\mu, & \beta_*\leq x<d.\\
\end{array}
\right.
\end{equation}
Otherwise we define
\begin{equation}\label{rzut_lc_ogolny}
P_yh_\alpha(x)=\left\{
\begin{array}{ll}
\dfrac{W_\alpha(y)-\hat{G}^{(k)}_n(y)}{1-W_\alpha(y)}\left(\dfrac{x-y}{\int\limits_0^yW_\alpha(x)dx}+1\right),&0\leq x<y,\\
&\\
\dfrac{W_\alpha(y)-\hat{G}^{(k)}_n(y)}{1-W_\alpha(y)},&y\leq x< d,
\end{array}
\right.
\end{equation}
with
\begin{equation}\nonumber
||P_yh_\alpha||^2_{W_\alpha}=\left(\frac{W_\alpha(y)-\hat{G}^{(k)}_n(y)}{(1-W_\alpha(y))\int\limits_0^y W_\alpha(x)dx}\right)^2\left(2y\int\limits_0^yW_\alpha(x)dx-2\int\limits_0^yxW_\alpha(x)dx-\left(\int\limits_0^yW_\alpha(x)dx\right)^2\right).
\end{equation}
Let $\mathcal{Z}$ denote the set of arguments $y\geq \beta_*$ satisfying the following condition
\begin{equation}\label{warunek_na_y_nasz}
\hat{G}^{(k)}_n(y)=W_\alpha(y)-(1-W_\alpha(y))\frac{\int\limits_0^yW_\alpha(x)dx\left(\int\limits_0^yW_\alpha(x)dx-\int\limits_0^y\hat{G}^{(k)}_n(x)dx\right)}{\left(2y-\int\limits_0^yW_\alpha(x)dx\right)\int\limits_0^yW_\alpha(x)dx-2\int\limits_0^yxW_\alpha(x)dx}.
\end{equation}
Then $\mathcal{Z}$ is nonempty and $P_y h_\alpha(x)=P_{y^*}h_\alpha(x)$ for unique $y^*=\arg\max\limits_{y\in\mathcal{Z}}||P_yh_\alpha||_{W_\alpha}$ , and we have the following bound
\begin{equation}\label{oszac_ogolne2}
\mathbb{E}\frac{R_n^{(k)}-\mu}{\sigma}\leq B_{\alpha,2}(y^*),
\end{equation}
where
\begin{eqnarray}
B_{\alpha,2}(y^*)&=&||P_{y^*}h_\alpha||_{W_\alpha}\nonumber\\
&=&\frac{W_\alpha(y^*)-\hat{G}^{(k)}_n(y^*)}{(1-W_\alpha(y^*))\int\limits_0^{y^*} W_\alpha(x)dx}\sqrt{2y^*\int\limits_0^{y^*}W_\alpha(x)dx-2\int\limits_0^{y^*}xW_\alpha(x)dx-\left(\int\limits_0^{y^*}W_\alpha(x)dx\right)^2}.\nonumber
\end{eqnarray}
The equality in \re{oszac_ogolne2} is attained for distributions $F\in$ IGFR($\alpha$) satisfying the following condition
\begin{equation}
F^{-1}(W_\alpha(x))=\left\{
\begin{array}{ll}
\dfrac{\sigma(W_\alpha(y^*)-\hat{G}^{(k)}_n(y^*))}{(1-W_\alpha(y^*))B_{\alpha,2}(y^*)}\left(\dfrac{x-y^*}{\int\limits_0^{y^*}W_\alpha(x)dx}+1\right) +\mu, &0\leq x<y^*,\\
 &\\
\dfrac{\sigma(W_\alpha(y^*)-\hat{G}^{(k)}_n(y^*))}{(1-W_\alpha(y^*)})B_{\alpha,2}(y^*) +\mu, & y^*\leq x<d.\\
\end{array}
\right.
\end{equation}

\end{proposition}

\begin{remark}
Proposition $\ref{Proposition_ogolne}$ describes general results for all possible parameters $\alpha$, with $k\geq 2$. Integrals that appear in expressions above strictly depend on $\alpha$ and since they have long analytic representations, we do not present them in the proposition, but gather them below,
\begin{eqnarray}
\int\limits_0^yW_\alpha(x)dx&=&\left\{
\begin{array}{ll}
-1+y+\mathrm{e}^{-y},&\alpha=0,\\
&\\
y+\dfrac{(1+\alpha y)^{1/\alpha+1}-1}{1+\alpha},&\alpha\neq 0,
\end{array}
\right.\nonumber\\
&&\nonumber\\
\int\limits_0^yxW_\alpha(x)dx&=&\left\{
\begin{array}{ll}
\frac{1}{2}y^2+\mathrm{e}^{-y}(1+y)-1,&\alpha=0,\\
&\\
\frac{1}{2}y^2+\dfrac{(1-\alpha y)^{1/\alpha+1}}{1+\alpha}y+\dfrac{(1-\alpha y)^{1/\alpha+2}-1}{(1+\alpha)(1+2\alpha)},&\alpha\neq 0,
\end{array}
\right.\nonumber\\
&&\nonumber\\
    \int_0^y\hat{G}_n^{(k)}(x)dx &=&\left\{
    \begin{array}{ll}
    y-\dfrac{1}{k}\sum\limits_{i=0}^n\left[1-{\rm e}^{-ky}\sum\limits_{j=0}^i\dfrac{(ky)^j}{j!}\right],&\alpha=0,\\
    &\\
    y-\sum\limits_{i=0}^n\dfrac{k^i}{(k+\alpha)^{i+1}}\left[1-\sum\limits_{j=0}^i(1-\alpha y)^{k/\alpha+1}\left[-\ln(1-\alpha y)^{1/\alpha}\right]^j\right.&\\
    \left.\cdot\dfrac{(k+\alpha)^j}{j!}\right],&\alpha\neq 0.
    \end{array}
    \right.\nonumber
\end{eqnarray}
Moreover, $\lambda_\alpha$ can be calculated from \re{lambda_alpha_ost} if $\lambda\neq 0$, and $\lambda_\alpha=\lambda_V$, further defined in \re{lambda_V}, if $\alpha=0$.
\end{remark}

\textsc{Proof of Proposition \ref{Proposition_ogolne}.}
Consider first case (i). Note that due to Lemma \ref{lemma_shape} in this setting of parameters, \re{h_alpha} is concave increasing and then decreasing. In this case the projection onto the cone of nondecreasing and concave functions is in fact the same as the projection onto the cone of nondecreasing functions. Indeed, the projection onto the family of nondecreasing functions coincides first with the original function on a subinterval of its concave increase, and then becomes constant. Therefore it is a nondecreasing concave function, and so it is the projection onto the cone of nondecreasing concave functions as well.
Hence
\begin{equation}\label{rzut_ogolny0}
P_\alpha h_\alpha(x)=\left\{
\begin{array}{ll}
h_\alpha(x),&0\leq x<\beta,\\
&\\
h_\alpha(\beta),& 0\leq \beta \leq d,
\end{array}
\right.
\end{equation}
for some $\beta\in(0,d)$. In order to determine parameter $\beta$ such that \re{rzut_ogolny0} is the projection of \re{h_alpha} onto \re{C_W}, we use condition \re{warunek_calka_rzutu} with
\begin{eqnarray}
\int\limits_0^dP_\alpha h_\alpha(x)w_\alpha(x)dx&=&\int\limits_0^\beta h_\alpha(x)w_\alpha(x)dx+h_\alpha(\beta)\int\limits_\beta^dw_\alpha(x)dx\nonumber\\
&=&\hat{G}^{(k)}_n(\beta)-W_\alpha(\beta)+(\hat{g}^{(k)}_n(\beta)-1)(1-W_\alpha(\beta)),\nonumber
\end{eqnarray}
which is equivalent to
\begin{equation}\label{warunek_beta1}
\hat{G}^{(k)}_1(\beta)=1-\hat{g}^{(k)}_1(\beta)(1-W_\alpha(\beta)),\quad 0\leq\beta<d,
\end{equation}
for $n=1$, and finally
$$
W_\alpha(\beta)=1-{\rm exp}\left(-\frac{1}{k(k-1)}\right),\quad 0\leq\beta<d,
$$
which allows to determine parameter \re{beta_star}. The bound \re{oszac_ogolne0} is determined by the norm of \re{rzut_ogolny0}.
\\

Note that for parameters $k, n, \alpha $, settled in case (ii), functions \re{h_alpha} satisfy conditions (A) and we directly use Proposition \ref{Proposition_GR}, presenting only a draft proof here.

Note that if $T_\alpha(b)<0$ and the set $\mathcal{Y}_\alpha$ is nonempty, then the projection of $h_\alpha$ is of the following form
\begin{equation}\nonumber
P_\alpha h_\alpha(x)=\left\{
\begin{array}{lll}
\hat{g}^{(k)}_n(y_*)-1+\lambda_\alpha(y_*)(x-y_*),&0\leq x<y_*,\\
&\\
\hat{g}^{(k)}_n(x)-1,&y_*\leq x<\beta_*,\\
&\\
\hat{g}^{(k)}_n(\beta_*)-1, &\beta_*\leq x<d.
\end{array}
\right.
\end{equation}
Otherwise, if these conditions are not satisfied, then the projection is just the linear increasing and constant function corresponding with formula \re{rzut_lc_ogolny}. Here the condition \re{warunek_na_y} turns out to be identical with \re{warunek_na_y_nasz}.

The bounds in (ii), according to \re{problem_Ph} can be determined by the norm of the projections $P_\alpha h_\alpha$, therefore we obtain \re{oszac_ogolne1} and finally \re{oszac_ogolne2} corresponding to $||P_{y^*}h_\alpha||_{W_\alpha}$, respectively.\\

The equality distributions for all considered cases can be determined using \re{warunek_rownosci} with $W=W_\alpha$ and $h=h_\alpha$ and appropriate projection functions $P_\alpha h_\alpha$. \hfill$\Box$\\


\begin{remark}
Note that for $n=1$ and $0<\alpha\leq 2k-1$, the construction of bounds implies that they are equal to the general ones, i.e. those derived without restricting to any special families of distributions $($see Raqab, $1997)$. The parameter \re{beta_star} is then the transformation of so called Moriguti point $\beta_0$ $($here equal to $0.3935)$, which defines the projection in the general case, according to the formula $\beta^*=W^{-1}_\alpha(\beta_0)$. Indeed, our condition $\re{warunek_beta1}$ for $\alpha=1$ is the same as the equation given by Raqab $(1997$, see formula $(2.4)$ for $n=2)$, with $\beta^*=W^{-1}_1(\beta_0)=U^{-1}(\beta_0)=0.3935$, matching parameter $u_1$ in Raqab's paper. For $\alpha=0$ we obtain transformed $\beta^*=W^{-1}_0(\beta_0)=V^{-1}(\beta_0)=0.5$, and the same value of the bound equal to $0.3451$ in both cases $($see Table $1$ below$)$.
\end{remark}

The particularly important cases of the distributions $F\in$ IGFR($\alpha$) are the distributions with the increasing density and increasing failure rate. The corresponding results are presented in the next subsections.

\subsection{Bounds for distributions with increasing density}

In case of the increasing density distributions (ID, for short) we fix $\alpha=1$ and $W_1(x)=U(x)=x$, $0\leq x<1$, as the standard uniform distribution, hence $\hat{g}^{(k)}_n(x)=g^{(k)}_n$ given by \re{g_nk} and
\begin{equation}\label{h_U}
  h_U(x)=h_1(x)=g^{(k)}_n(x)-1.
\end{equation}

Case $k=1$, when \re{h_U} is convex increasing (see Rychlik (2001, p.136) and Lemma \ref{lemma_shape}), has already been considered in Proposition \ref{Proposition_k1}. Here we consider other cases, when \re{h_U} is concave increasing on $(0,1-\exp(-\frac{1}{k-1}))$ and decreasing on $(1-\exp(-\frac{1}{k-1}),1)$ for $k\geq 2$, $n=1$, and  when \re{h_U} satisfy conditions (A) for $k\geq 2$, $n\geq 2$ with
\begin{eqnarray}
b&=&b_U=\left\{
\begin{array}{ll}
  1-\exp(-(n-1)), & k=2, \\
  &\\
  1-\exp\left(-\dfrac{2kn-3n-\sqrt{n^2+4n(k-1)(k-2)}}{2(k-1)(k-2)}\right), & k\geq 3,
\end{array}
\right.\nonumber\\
c&=&c_U=1-\exp\left(-\frac{n}{k-1}\right)\nonumber.
\end{eqnarray}
Note that functions defined in \re{T}-\re{Z} take the following form in the ID case of distributions
\begin{eqnarray}
T_U(\beta)&=&(1-\beta)\left[-\frac{1}{k}\sum\limits_{i=0}^{n-1}g^{(k)}_i(\beta)+\left(1-\frac{1}{k}\right)g^{(k)}_n(\beta)\right],\label{T_U}\\
\lambda_U(y)&=& \frac{3}{y^3}\left\{ \sum\limits_{i=0}^n\frac{k^i}{(k+1)^{i+1}}G^{(k+1)}_i(y)+\frac{\alpha^2}{2}g^{(k)}_n(y)-y \right\}\nonumber\\
&=&\frac{3}{y^3}\left\{1-\alpha+\frac{\alpha^2}{2}g^{(k)}_n(y)-\left(\frac{k}{k+1}\right)^{n+1}-\sum\limits_{i=0}^n\frac{k^i}{(k+1)^{i+1}}\left(\sum\limits_{j=0}^i\frac{(k+1)^j}{(k+2)^{j+1}}g^{(k+2)}_j(y)\right) \right\},\nonumber\\
Y_U(y) &=& \lambda_U(y) -\frac{k^{n+1}}{(k-1)^n}\left(g^{(k-1)}_{n-1}(y)-g^{(k-1)}_{n}(y)\right),\nonumber \\
Z_U(y)&=&1-yg^{(k)}_n(y)-\sum\limits_{i=0}^n\frac{k^i}{(k+1)^{i+1}}g^{(k+1)}_i(y)\nonumber\\
&&+\frac{3}{2y}\left\{1-\alpha+\frac{\alpha^2}{2}g^{(k)}_n(y)-\left(\frac{k}{k+1}\right)^{n+1}-\sum\limits_{i=0}^n\frac{k^i}{(k+1)^{i+1}}\left(\sum\limits_{j=0}^i\frac{(k+1)^j}{(k+2)^{j+1}}g^{(k+2)}_j(y)\right) \right\},\nonumber
\end{eqnarray}
for $0\leq \beta,y<1$.\\

The results below follow from the Proposition \ref{Proposition_ogolne} for $\alpha=1$.
\begin{corollary}\label{Proposition_ID}
Let $F$ be the increasing density distribution function and fix $k\geq 2$. \\
(i) For $n=1$ let $\beta^*$ satisfy the following condition
\begin{equation}\nonumber
G^{(k)}_1(\beta)=g^{(k)}_1(\beta)(\beta-1)+1,\quad 0\leq \beta\leq 1.
\end{equation}
 Then
 \begin{equation}\label{oszac_ID_0}
 \mathbb{E}\frac{R_1^{(k)}-\mu}{\sigma}\leq B_{U,0}(\beta^*),
 \end{equation}
 where
 \begin{equation}\nonumber
 B_{U,0}(\beta)=\frac{2k^4}{(2k-1)^3}G^{(2k-1)}_2(\beta)-2G^{(k)}_1(\beta)+\beta+(1-\beta)(g^{(k)}_1(\beta)-1)^2.
 \end{equation}
 The equality in \re{oszac_ID_0} is attained for the following distribution function
 \begin{equation}\nonumber
 F(x)=\left\{
 \begin{array}{ll}
 0,&x<\mu-\frac{\sigma}{B_{U,0}(\beta^*)},\\
 &\\
 \left( g^{(k)}_1\right)^{-1}(B_{U,0}(\beta^*)\frac{x-\mu}{\sigma}+1),& \mu-\frac{\sigma}{B_{U,0}(\beta^*)}\leq x<\frac{g^{(k)}_1(\beta^*)-1}{B_{U,0}(\beta^*)}\sigma+\mu,\\
&\\
1,&x\geq \frac{g^{(k)}_1(\beta)-1}{B_{U,0}(\beta_*)}\sigma+\mu.
 \end{array}
 \right.
 \end{equation}

(ii) Let now $n\geq 2$ and $\beta_U$ denote the unique zero of \re{T_U} in the interval $(0,c_U)$. If $T_U(b_U)<0$ and the set $\mathcal{Y_U}=\{y\in(b_U,\beta_U): Y_U(y)\geq 0, Z_U(y)=0\}$ is nonempty, then
\begin{equation}\nonumber
P_Uh_U(x)=\left\{
\begin{array}{lll}
g^{(k)}_n(y_U)-1+\lambda_U(y_U)(x-y_U),&0\leq x<y_U,\\
&\\
g^{(k)}_n(x)-1,&y_U\leq x<\beta_U,\\
&\\
g^{(k)}_n(\beta_U)-1, &\beta_U\leq x<1,
\end{array}
\right.
\end{equation}
where $y_U=\inf\{y\in \mathcal{Y_U}\}$ and we have the following bound
\begin{equation}\label{oszac_ID_1}
\mathbb{E}\frac{R_n^{(k)}-\mu}{\sigma}\leq B_{U,1}(y_U,\beta_U),
\end{equation}
where
\begin{eqnarray}
B_{U,1}^2(y,\beta)&=&y\left[\left(g^{(k)}_n(y)\right)^2-\left(y\lambda_U(y)+2\right)g^{(k)}_n(y)+\frac{1}{3}y\lambda_U(y)\left(y\lambda_U(y)+3\right)\right]\nonumber\\
&&+\beta-2\left[G^{(k)}_n(\beta)-G^{(k)}_n(y)\right]+(1-\beta)\left(g^{(k)}_n(\beta)-1\right)^2\nonumber\\
&&+\frac{k^{2(n+1)}(2n)!}{(n!)^2(2k-1)^{2n+1}}\left[G^{(2k-1)}_{2n}(\beta)-G^{(2k-1)}_{2n}(y)\right].\nonumber
\end{eqnarray}
The equality in \re{oszac_ID_1} holds for the following distribution function
\begin{equation}\nonumber
F(x)=\left\{
\begin{array}{ll}
0,& x<\frac{g^{(k)}_n(y)-1-y\lambda_U(y)}{B}\sigma+\mu,\\
&\\
y+\dfrac{1}{\lambda_U(y)}\left(B\frac{x-\mu}{\sigma}-g^{(k)}_n(y)+1\right),& \frac{g^{(k)}_n(y)-1-y\lambda_U(y)}{B}\sigma+\mu\leq x\\
&<\frac{g^{(k)}_n(y)-1}{B}\sigma+\mu,\\
&\\
\left(g^{(k)}_n\right)^{-1}(1+\frac{x-\mu}{\sigma}B),&\frac{g^{(k)}_n(y)-1}{B}\sigma+\mu\leq x<\frac{g^{(k)}_n(\beta)-1}{B}\sigma+\mu,\\
&\\
1,&x\geq \frac{g^{(k)}_n(\beta)-1}{B}\sigma+\mu,
\end{array}
\right.
\end{equation}
with $B=B_{U,1}(y_U,\beta_U)$, $y=y_U$ and $\beta=\beta_U$.
Otherwise let
\begin{equation}\nonumber
P_yh_U(x)=\left\{
\begin{array}{ll}
\dfrac{2(y-G^{(k)}_n(y))}{y^2(1-y)}\cdot x-\dfrac{(2-y)(y-G^{(k)}_n(y))}{y(1-y)},  & 0\leq x<y, \\
&\\
\dfrac{y-G^{(k)}_n(y)}{1-y}, & y\leq x<1,
\end{array}
\right.
\end{equation}
with
\begin{equation}\nonumber
||P_yh_U||^2_U=\frac{4-3y}{3y}\left(\frac{y-G^{(k)}_n(y)}{1-y}\right)^2.
\end{equation}
Let $\mathcal{Z}_U$ denote the set of $y\geq \beta_U$ being the solution to the following equation
\begin{equation}\nonumber
G^{(k)}_n(y)=y-\frac{6(1-y)}{y(4-3y)}\left(\frac{y^2}{2}-y+\sum\limits_{i=1}^n\frac{k^i}{(k+1)^{i+1}}G^{(k+1)}_i(y)\right).
\end{equation}
Then $\mathcal{Z}_U$ is nonempty and $P_yh_U(x)=P_{y_U^*}h_U(x)$ for unique $y_U^*=\arg\max\limits_{y\in\mathcal{Z}_U}||P_yh_U||_U$ and we have the following bound
\begin{equation}\label{oszac_ID_2}
\mathbb{E}\frac{R_n^{(k)}-\mu}{\sigma}\leq B_{U,2}(y_U^*),
\end{equation}
with
\begin{equation}\nonumber
B_{U,2}(y_U^*)=\frac{y_U^*-G^{(k)}_n(y_U^*)}{1-y_U^*}\sqrt{\frac{4-3y_U^*}{3y_U^*}}.
\end{equation}
The equality in \re{oszac_ID_2} holds for the following distribution function
\begin{equation}\nonumber
F(x)=\left\{
\begin{array}{ll}
0,&x<\frac{\sigma}{B}A(1-\frac{y_U^*}{C})+\mu,\\
&\\
C(\frac{x-\mu}{\sigma}\frac{B}{A}-1)+y_U^*,&\frac{\sigma}{B}A(1-\frac{y_U^*}{C})+\mu\leq x<A\frac{\sigma}{B}+\mu,\\
&\\
1,&x\geq A\frac{\sigma}{B}+\mu,
\end{array}
\right.
\end{equation}
with
\begin{eqnarray}
A&=&\frac{y_U^*-G^{(k)}_n(y_U^*)}{1-y_U^*},\nonumber\\
B&=&B_{U,2}(y_U^*),\nonumber\\
C&=&\frac{1}{2}(y_U^*)^2.\nonumber
\end{eqnarray}

\end{corollary}

\subsection{Distributions with increasing failure rate}
Let us now consider distibutions with the increasing failure rate (IFR, for short), i.e. $\alpha=0$ and $W_0(x)=V(x)=1-\textrm{e}^{-x}$, $0\leq x<\infty$, which is the standard exponential distribution. Therefore we have
\begin{eqnarray}
  \hat{g}^{(k)}_n(x)&=&g^{(k)}_n(V(x))=\frac{k^{n+1}}{n!}x^n\textrm{e}^{-x(k-1)},\nonumber\\
  \hat{G}^{(k)}_n(x)&=&1-\textrm{e}^{-kx}\sum\limits_{i=0}^n\frac{(xk)^i}{i!}=1-\sum\limits_{i=1}^n\frac{k^i}{(k+1)^{i+1}}\hat{g}^{(k+1)}_i(x),\nonumber\\
  h_V(x)&=&h_0(x)=\hat{g}^{(k)}_n(x)-1.\label{h_V}
\end{eqnarray}
The case of the first records ($k=1$) with $n\geq 1$ for the increasing failure rate distributions was presented in Proposition \ref{Proposition_k1} (ii). Therefore below we consider only $k\geq 2$ with $n\geq 1$. Here
\begin{eqnarray}
b&=&b_V=\frac{n-\sqrt{n}}{k-1},\nonumber\\
c&=&c_V=\frac{n}{k-1},\nonumber
\end{eqnarray}
and the respective functions \re{T}-\re{Z} are following
\begin{eqnarray}
T_V(\beta)&=&-\sum\limits_{i=0}^{n-1}\frac{k^i}{(k+1)^i+1}\hat{g}^{(k+1)}_i(\beta)+(k-1)\frac{k^n}{(k+1)^n}\hat{g}^{(k+1)}_n(\beta),\label{T_V}\\
\lambda_V(y)&=& \frac{1}{y^2+2(1-y-\textrm{e}^{-y})}\left\{\left(y-\frac{n+1}{k}\right)\left[\sum\limits_{i=0}^n\frac{k^i}{(k+1)^{i+1}}\hat{g}^{(k+1)}_i(y)-1\right]\right.\label{lambda_V}\\
&&\left.-\frac{(n+1)k^n}{(k+1)^{n+2}}\hat{g}^{(k+1)}_{n+1}(y)-(1-y-\textrm{e}^{-y})\hat{g}^{(k)}_n(y)\right\},\nonumber\\
Y_V(y)&=&\lambda_V(y)-[k\hat{g}^{(k)}_{n-1}(y)-(k-1)\hat{g}^{(k)}_n(y)],\nonumber\\
Z_V(y)&=&1-\sum\limits_{i=0}^n\frac{k^i}{(k+1)^{i+1}}\hat{g}^{(k+1)}_i(y)-(1-\textrm{e}^{-y})\hat{g}^{(k)}_n(y)-\lambda_V(y)(1-y-\textrm{e}^{-y}).\nonumber
\end{eqnarray}
The results below are the straightforward implication of Proposition \ref{Proposition_ogolne}, therefore the proof of the corollary below is immediate and will not be presented here.
\begin{corollary}
Let $F$ be the increasing failure rate cdf and fix $k\geq 2$.

 (i) Let $n=1$ and $\beta^*$ satisfy condition
 \begin{equation}\nonumber
 \hat{G}^{(k)}_n(\beta)=1-\textrm{e}^{-\beta}\hat{g}^{(k)}_1(\beta),\quad 0\leq \beta<\infty.
 \end{equation}
 Then we have the following bound
 \begin{equation}\nonumber
 \mathbb{E}\frac{R_1^{(k)}-\mu}{\sigma}\leq B_{V,0}(\beta^*),
 \end{equation}
 where
 \begin{equation}\nonumber
 B_{V,0}^2(\beta)=\frac{2k^4}{(2k-1)^3}\hat{G}^{(2k-1)}_2(\beta)-2\hat{G}^{(k)}_1(\beta)+1-\textrm{e}^{-\beta}[1-(\hat{g}^{(k)}_1(\beta)-1)^2].
 \end{equation}
The equality above is attained for the following distribution function
\begin{equation}\nonumber
F(x)=\left\{
\begin{array}{ll}
0,& x<\mu-\frac{\sigma}{ B_{V,0}(\beta^*)},\\
&\\
1-\exp\left(-\left(\hat{g}^{(k)}_1\right)^{-1}\left( B_{V,0}(\beta^*)\dfrac{x-\mu}{\sigma}+1\right)\right),&\mu-\dfrac{\sigma}{ B_{V,0}(\beta^*)}\leq x\\
&<\dfrac{\sigma}{ B_{V,0}(\beta^*)}(\hat{g}^{(k)}_1(\beta^*)-1)+\mu,\\
&\\
1,&x\geq \dfrac{\sigma}{ B_{V,0}(\beta^*)}(\hat{g}^{(k)}_1(\beta^*)-1)+\mu.
\end{array}
\right.
\end{equation}
(ii) Let $n\geq 2$ and let $\beta_V$ denote the unique zero of \re{T_V}. If $T_V(b_V)<0$ and the set $\mathcal{Y_V}=\{y\in(b_V,\beta_V): Y_V(y)\geq 0, Z_V(y)=0\}$ is nonempty, then
\begin{equation}\nonumber
P_Vh_V(x)=\left\{
\begin{array}{lll}
\hat{g}^{(k)}_n(y_V)-1+\lambda_V(y_V)(x-y_V),&0\leq x<y_V,\\
&\\
\hat{g}^{(k)}_n(x)-1,&y_V\leq x<\beta_V,\\
&\\
\hat{g}^{(k)}_n(\beta_V)-1, &x\geq\beta_V,
\end{array}
\right.
\end{equation}
where $y_V=\inf\{y\in \mathcal{Y_V}\}$ and we have the following bound
\begin{equation}\label{oszac_IFR_1}
\mathbb{E}\frac{R_n^{(k)}-\mu}{\sigma}\leq B_{V,1}(y_V,\beta_V),
\end{equation}
where
\begin{eqnarray}
B_{V,1}^2(y,\beta)&=&(1-\mathrm{e}^{-y})\left(\hat{g}^{(k)}_n(y)\right)^2-2\hat{g}^{(k)}_n(y)[1-\mathrm{e}^{-y}-\lambda_V(y)(1-y-\mathrm{e}^{-y})]\nonumber\\
&&+\lambda_V^2(y)[y^2+2(1-y-\mathrm{e}^{-y})]-2\lambda_V(y)(1-y-\mathrm{e}^{-y})+1+\mathrm{e}^{-\beta}[(\hat{g}^{(k)}_n(\beta)-1)^2 -1]\nonumber\\
&&+\dfrac{k^{2(n+1)}(2n)!}{(n!)^2(2k-1)^{2n+1}}\left(\hat{G}^{(2k-1)}_{2n}(\beta)-\hat{G}^{(2k-1)}_{2n}(y)\right)-2\left(\hat{G}^{(k)}_n(\beta)-\hat{G}^{(k)}_n(y)\right).\nonumber
\end{eqnarray}
The equality in \re{oszac_IFR_1} holds for the following distribution function
\begin{equation}\nonumber
F(x)=\left\{
\begin{array}{ll}
0,&x<\frac{\sigma}{B(y,\beta)}(\hat{g}^{(k)}_n(y)-1-y\lambda_V(y))+\mu,\\
&\\
V\left(\frac{1}{\lambda_V(y)}(B\frac{x-\mu}{\sigma}-\hat{g}^{(k)}_n(y)+1)+y\right),&\frac{\sigma}{B}(\hat{g}^{(k)}_n(y)-1-y\lambda_V(y))+\mu\leq x\\
&< \frac{\sigma}{B}(\hat{g}^{(k)}_n(y)-1)+\mu,\\
&\\
V\left(\left(\hat{g}^{(k)}_n\right)^{-1}(B\frac{x-\mu}{\sigma}+1)\right),&\frac{\sigma}{B}(\hat{g}^{(k)}_n(y)-1)+\mu \leq x\\
&< \frac{\sigma}{B}(\hat{g}^{(k)}_n(\beta)-1)+\mu ,\\
&\\
1,&x\geq\frac{\sigma}{B}(\hat{g}^{(k)}_n(\beta)-1)+\mu ,
\end{array}
\right.
\end{equation}
with $B=B_{V,1}(y_V,\beta_V)$, $y=y_V$ and $\beta=\beta_V$.
Otherwise let
\begin{equation}\nonumber
P_yh_V(x)=\left\{
\begin{array}{ll}
\dfrac{1-\mathrm{e}^{-y}-\hat{G}^{(k)}_n(y)}{\mathrm{e}^{-y}}\left(\dfrac{x-y}{\mathrm{e}^{-y}+y-1}+1\right),  & 0\leq x<y, \\
&\\
\dfrac{1-\mathrm{e}^{-y}-\hat{G}^{(k)}_n(y)}{\mathrm{e}^{-y}}, & x\geq y,
\end{array}
\right.
\end{equation}
with
\begin{equation}\nonumber
||P_yh_V||^2_V=\left(\frac{1-\mathrm{e}^{-y}-\hat{G}^{(k)}_n(y)}{\mathrm{e}^{-y}+y-1}\right)^2(\mathrm{e}^{2y}-1-2y\mathrm{e}^y).
\end{equation}
Let $\mathcal{Z}_V$ denote the set of $y\geq \beta_V$ which satisfy the following equation
\begin{eqnarray}
\hat{G}^{(k)}_n(y)&=&1-\mathrm{e}^{-y}+\frac{\mathrm{e}^{-y}(1-y-\mathrm{e}^{-y})}{1-\mathrm{e}^{-y}(\mathrm{e}^{-y}+2y)}\left[\frac{n+1}{k}-1+\mathrm{e}^{-y}\right.\nonumber\\
&&\left.+ \left(y-\frac{n+1}{k}\right)\sum\limits_{i=0}^n\frac{k^i}{(k+1)^{i+1}}\hat{g}^{(k+1)}_i(y)-\frac{(n+1)k^n}{(k+1)^{n+2}}\hat{g}^{(k+1)}_{n+1}(y) \right].\nonumber
\end{eqnarray}
Then $\mathcal{Z}_V$ is nonempty and $P_yh_V(x)=P_{y_V^*}h_V(x)$ for unique $y_V^*=\arg\max\limits_{y\in\mathcal{Z}_V}||P_yh_V||_V$ and we have the following bound
\begin{equation}\label{oszac_IFR_2}
\mathbb{E}\frac{R_n^{(k)}-\mu}{\sigma}\leq B_{V,2}(y_V^*),
\end{equation}
where
\begin{equation}\nonumber
B_{V,2}(y)=||P_{y}h_V||_V=\frac{1-\mathrm{e}^{-y}-\hat{G}^{(k)}_n(y)}{\mathrm{e}^{-y}+y-1}\sqrt{\mathrm{e}^{2y}-1-2y\mathrm{e}^{y}}.
\end{equation}
The equality in \re{oszac_IFR_2} holds for the following distribution function
\begin{equation}\nonumber
F(x)=\left\{
\begin{array}{ll}
0,&x<\frac{\sigma}{B}A(1-\frac{y_V^*}{C})+\mu,\\
&\\
V\left(C(\frac{x-\mu}{\sigma}\frac{B}{A}-1)+y_V^*\right),&\frac{\sigma}{B}A(1-\frac{y_V^*}{C})+\mu\leq x<A\frac{\sigma}{B}+\mu,\\
&\\
1,&x\geq A\frac{\sigma}{B}+\mu,
\end{array}
\right.
\end{equation}
with
\begin{eqnarray}
A&=&\frac{1-\mathrm{e}^{-y_V^*}-\hat{G}^{(k)}_n(y_V^*)}{\mathrm{e}^{-y_V^*}},\nonumber\\
B&=&B_{V,2}(y_V^*),\nonumber\\
C&=&\mathrm{e}^{-y_V^*}+y_V^*-1.\nonumber
\end{eqnarray}
\end{corollary}

\subsection{Numerical calculations}
We illustrate the obtained results with the numerical calculations of the bounds on \re{problem} for fixed values of parameters $n$, $k$ and $\alpha$.

Table 1 below contains the numerical values of the upper bounds on the $n$th values of the second records ($k=2$), $n=1,\ldots,9$, in two particular cases of the IGFR distributions: $F\in$ ID and $F\in$ IFR. The column $y$ presents the values of $y_U, y_V$ in case of the l-h-c type of the projection (inclined font of $y$) and $y^*_U,y^*_V$, if the projection has the shape l-c. For $n=1$ column $y$ presents $\beta^*$ (see Corollary 1 and 2, cases (i)) for particular shape of projection (h-c, say). As it was mentioned before (see Remark 2 above), for $n=1$ we obtain the same value $y=\beta^*=0.3935$ as Raqab (1997) obtained for $u_1$ in case of ID distribution, and for the IFR family the parameter is transformed, $y=\beta^*=V^{-1}(0.3935)=0.5$. Moreover, Raqab's numerical value of the bound $B_{2,2}(2)=0.34507$ is the same as our bounds for both ID and IFR families.

The only case among those considered when the projection shape is l-h-c, is the case $F\in$ IFR, $n=2$, where $\beta_V=1.3660$. Moreover, for $F\in$ ID and $n=6,\ldots,9$, the value $y=1.000$ is only the approximation of the change point which approaches to 1, from which the projection becomes constant.

Note that the bounds increase while $n$ increases, and so do the parameters $y$. Moreover, the bounds in the IFR case are greater than the bounds in the ID case, which is consistent with the general dependencies between those two families of distributions.
\begin{table}[ht]
\caption{Upper bounds on the expectations of standardized $n$th values of the second records from restricted families}
\begin{tabular}{c||c|c||c|c}
&\multicolumn{2}{|c||}{$F\in $ID} & \multicolumn{2}{c}{$F\in$ IFR}\\\hline
$n$&$y$&bound&$y$&bound\\\hline
1&0.3935 &0.3451& 0.5000& 0.3451\\
2&0.7954&0.7270&\textit{1.1433}&0.7350\\
3&0.9696&1.0485&2.3791&1.1321\\
4&0.9972&1.2759&3.6664&1.5600\\
5&0.9998&1.4280&5.1766&2.0214\\
6&1.0000&1.5293&6.9094&2.5059\\
7&1.0000&1.5969&8.8328&3.0013\\
8&1.0000&1.6417&10.8987&3.5002\\
9&1.0000&1.6713&13.0648&4.0000\\
\end{tabular}
\end{table}

Table 2 presents the optimal bounds on the expectations of the fifth values ($n=5$) of the $k$th records, $k=1,2,\ldots,10$, based on the distributions with the increasing density (ID) and increasing failure rate (IFR) respectively. For $k=1$ Proposition \ref{Proposition_k1} (ii) and (iii) was usefull.
For $k\geq 2$ all the cases of the bounds were determined based on the second kind of the projection shape, l-c (Corollary 1 and 2, cases (ii)). Note that the bounds decrease along with the increase of the parameter $k$, and the same concerns the parameters $y^*_U$, $y^*_V$, determining the points in which the projection breaks from the linear increasing into the constant function. It is not possible deliver parameter $y^*_W$, $W=U,V$ for $n=1$, since the shapes of \re{h_U} and \re{h_V} are then linear increasing.

\begin{table}[ht]
\caption{Upper bounds on the expectations of standardized fifth values of the $k$th records from restricted families}
\begin{tabular}{c||c|c||c|c}
&\multicolumn{2}{|c||}{$F\in $ID} & \multicolumn{2}{c}{$F\in$ IFR}\\\hline
$k$&$y^*_U$&bound&$y^*_V$&bound\\\hline
1 &-&1.6779&-&5\\
2&0.9998&1.4279&5.1766&2.0214\\
3&0.9328&1.1209&2.1472&1.2296\\
4&0.7672&0.8875&1.3001&0.9209\\
5&0.6226&0.7389&0.9087&0.7544\\
6&0.5137&0.6393&0.6870&0.6482\\
7&0.4322&0.5678&0.5460&0.5736\\
8&0.3701&0.5137&0.4493&0.5179\\
9&0.3217&0.4711&0.3793&0.4743\\
10&0.2832&0.4366&0.3266&0.4391\\
\end{tabular}
\end{table}
\section*{Acknowledgements}
The author is greatly indebted to anonymous reviewers for many valuable comments which helped in the preparation of the final version of the paper. The research was supported by the Polish National Science Center Grant no. 2015/19/B/ST1/03100.

\end{document}